\documentclass[reqno, 12pt]{amsart}
\usepackage{hyperref}
\usepackage{amsmath, amsthm, amssymb, mathtools, graphicx, paralist, color}
\usepackage{indentfirst}
\usepackage{appendix}
\usepackage{listings}
\usepackage{hyphenat}
\usepackage{calrsfs}
\usepackage{enumerate}
\usepackage{tikz-cd}
\usepackage{mathrsfs}
\usepackage{comment}
\hypersetup{
  hidelinks,
  pdftitle={Square-root cancellation in three hyperbolic counting problems},
  pdfauthor={Marios Voskou}
}

\newtheorem{corollary}{Corollary}
\newtheorem{theorem}{Theorem}
\newtheorem{lemma}{Lemma}
\newtheorem{proposition}{Proposition}
\theoremstyle{remark}
\newtheorem{remark}{Remark}
\usepackage[utf8]{inputenc}

\def\R{\ensuremath\mathbb{R}}
\def \H{\ensuremath \mathbb{H}^2}
\def \Hh{\ensuremath \mathbb{H}^3 }

\def \psl  {{\hbox{PSL}_2( {\mathbb R})} }

\newcommand{\pslr}{\mathrm{PSL}_{2}(\mathbb{R})}

\newcommand{\pslc}{\mathrm{PSL}_{2}(\mathbb{C})}
\DeclareDocumentCommand{\pd}{O{} O{} m}{\frac{\partial^{#1}#2}{\partial{#3}^{#1}}}
\usepackage{bbm}
\usepackage[
backend=biber,
style=alphabetic,
sorting=nyt,
maxbibnames=99,
isbn=false,
doi=false
]{biblatex}
\title[Square-root cancellation in three hyperbolic counting problems]{Square-root cancellation in three hyperbolic counting problems via maximal estimates}
\author{Marios Voskou }
\address{Alfr\'{e}d R\'{e}nyi Institute of Mathematics, POB 127, Budapest H-1364, Hungary}
\email{voskou.marios@renyi.hu}
\date{\today}

\begin{document}
\begin{abstract} 
We study three counting problems for lattice orbits in hyperbolic space: counting orbit points within a prescribed distance of a fixed closed geodesic or a totally geodesic hyperplane, and counting orbit points above a prescribed height at a cusp. For almost every base point, we obtain error bounds of square-root size up to logarithmic factors. For the geodesic and hyperplane counts, these bounds match the conjectured exponents and improve earlier results due to Good and Laaksonen. For these two problems, we adapt the method of maximal estimates developed by Chatzakos and Dimakis for hyperbolic lattice counting, formulating it in both a discrete and a continuous framework. For the hyperplane problem, we also obtain an asymptotic formula for the spatial second moment of the error term, involving an almost periodic function. For the cusp problem, the square-root bound can be deduced from work of Yu, which incorporates Schmidt's dyadic maximal argument. We present a proof, together with an $\Omega$-result supporting the optimality of the exponent $1/2$ for at least one cusp of any group that is not essentially cuspidal. \end{abstract}
\maketitle
\section{Introduction}

For a fixed cofinite Fuchsian group $\Gamma \subset \psl$ and points $z,w \in \mathbb{H}^2$, where $\mathbb{H}^2$ is the hyperbolic upper half-plane, the associated hyperbolic lattice counting problem is concerned with estimating the number \begin{equation*} \label{clcountingf} N(X;z,w):= \# \left\{ \gamma \in \Gamma \vert  2\cosh{\rho(\gamma z,w) } \leq X \right\}.\end{equation*}
Here, $\rho$ is the hyperbolic distance function, induced from the Poincar\'e metric $ds(z)^2=y^{-2}\left( dx^2+dy^2 \right)$, where $z=x+iy$. The asymptotic behaviour of the function $N(X;z,w)$ is a classical problem in analytic number theory, initially studied by Selberg \cite{selberg} (see also the work of G\"{u}nther \cite{Gunther1980} and Good \cite{Good1983}). In particular, Selberg proved that
$$N(X;z,w)=\sum_{1/2<s_j \leq 1}a_j(z,w)X^{s_j}
+E(X; z,w),$$
where the sum is over the finitely many small eigenvalues $\lambda_j=s_j(1-s_j)$ of the Laplacian over $L^2(\Gamma \backslash \H)$, and $E(X;z,w)=O_{z,w}(X^{2/3})$. For explicit expressions for the coefficients $a_j(z,w)$, 
see \cite[equation (1.1)]{ChatzakosPetridis+2016+981+1003}. Furthermore, it is conjectured that, for every $\epsilon>0$, we have $E(X; z,w)=O_{z,w,\epsilon}(X^{1/2+\epsilon})$ (see \cite[p.81]{philipsr}). Until very recently, no pointwise progress towards this conjecture had been made, for any choice of $z,w$ and $\Gamma$. We now have the following essentially optimal improvement for almost all $z,w$, and for all cocompact groups $\Gamma$.
\begin{theorem}[Chatzakos--Dimakis \cite{chatzakos2026squarerootcancellationhyperbolic}] \label{almostall} Let $\Gamma$ be a cocompact Fuchsian group and fix a point $w \in \mathbb{H}^2$. Then, there exists a set $B_w \subset \Gamma \backslash \mathbb{H}^2$ of full hyperbolic measure, such that for every $z \in B_w$ and for every $\epsilon >0 $, we have
$$
E(X;z,w) = O_{\Gamma, w,z,\epsilon}  \left(X^{1/2} (\log X)^{3/2} (\log \log X)^{1/2+\epsilon}\right)
$$
as $X \to \infty$.
\end{theorem}

\begin{remark}
We note that while their result gives us an essentially optimal bound for almost all pairs of $z,w$, it does not provide any explicit pairs of $z,w$. For a logarithmic improvement for a particular choice of $z,w,\Gamma$, see \cite{chatzakos2025hyperboliccircleproblemheegner}.
\end{remark}
Their proof provides a new way to study sums of the form $\sum_{j}c_j(R)\phi_j(z)$ for $\mu$-almost every $z \in \mathcal{X}$, where $(\mathcal{X},\mu)$ is a measure space and $(\phi_j)_j$ is an orthonormal sequence of functions. Their method is particularly useful in situations where we have good mean-square bounds for the coefficients and their $R$-derivatives. Such sums appear extensively in spectral theory, and in particular in its use in studying counting problems. In this paper, we apply their methods to investigate two such problems, providing improvements on previous results, which we expect to be essentially optimal. For further applications, see also \cite{chatzakos2026pointwiseimprovementerrorterm} and \cite{dimakis2026offdiagonalboundsspectralprojector}. We also investigate a third counting problem, using maximal estimates in a different way.
\subsection{The point-geodesic problem}
For the first problem, instead of considering the distances of the images of a point $z \in \mathbb{H}^2$ from a point $w$, we will consider their distances from an oriented closed geodesic $l \subset \Gamma \backslash \mathbb{H}^2$, which lifts to the oriented line $\tilde{l} \subset \mathbb{H}^2$. Let $h$ be a generator of $\hbox{Stab}(\tilde{l})=:H$, so that $l=H\backslash \tilde{l}$. We study the number
$$N(X;z,l):= \# \left\{ \gamma \in H \backslash\Gamma \vert  \cosh{\rho( \gamma z, \tilde{l} ) } \leq X \right\}.$$
We can think of this as studying the distribution of cosets $H\gamma$ in $H \backslash \Gamma$, with respect to an appropriate parameterisation, where $H$ is any hyperbolic cyclic subgroup. Furthermore, this is also equivalent to the lattice counting problem restricted to the conjugacy class of $h$ (see \cite[\S 1]{ChatzakosPetridis+2016+981+1003}).
As in the classical hyperbolic lattice counting problem, we can write
$$N(X;z,l)=\sum_{1/2<s_j \leq 1}a_j(z,l)X^{s_j}+E(X;z,l),$$
where the coefficients $a_j(z,l)$ are as in equation (\ref{ajzldef}).
In \cite{Good1983}, Good proved that the error term satisfies $E(X;z,l)=O_{z,l}(X^{2/3})$. This was later reproved by Chatzakos--Petridis in  \cite{ChatzakosPetridis+2016+981+1003}, using a cleaner method, based on a weaker result of Huber \cite{Huber1998}. We prove the following improvement to their theorem, which reaches the expected square-root scale, up to logarithmic factors:
\begin{theorem} \label{pointline}
    Let $\Gamma$ be a cofinite Fuchsian group, and fix a primitive oriented closed geodesic $l$ in $\Gamma \backslash \H$. Then, there exists a set $B_l \subset \Gamma \backslash \mathbb{H}^2$ of full hyperbolic measure such that, for every $z \in B_l$, we have
$$
E(X;z,l) = O_{\Gamma, z,l}  \left(X^{1/2} (\log X)^{3/2} (\log \log X)^{1/2}\log\log\log X\right)
$$
as $X \to \infty$.
\end{theorem} 
\begin{remark}

After subtracting the contribution of the eigenvalue $1/4$, Chatzakos \cite[Theorem 1.13]{cha19} proves that, for every cocompact group and fixed closed geodesic, there is a point on that geodesic at which the error is $\Omega_+(X^{1/2})$. Thus, the exponent $1/2$ cannot be improved in a bound valid at every point.
    \end{remark}
    \begin{remark}As is the case for \cite{chatzakos2026squarerootcancellationhyperbolic}, we do not know whether the exponents of $\log{X}$ and $\log\log{X}$ are optimal. The exponent for $\log\log\log{X}$ is not. In fact, better exponents can be proved directly from our method. It would be interesting to investigate whether the exponent of $\log{X}$, after subtracting the potential contribution of $\lambda=1/4$, can be reduced below $1$. In that case, when the contribution of $\lambda=1/4$ is nonzero, it would dominate the rest of the error term. This follows from the fact that this contribution is $\sim c X^{1/2}\log{X}$, for some real number $c$.
\end{remark}

\subsection{The point-hyperplane problem}
In the second problem, we consider the distance of the images of a point in $\Hh $ from a fixed hyperplane. This can be viewed as a three-dimensional analogue of the first problem.
Let $\Gamma$ be a cocompact discrete subgroup of $\mathrm{PSL}_2(\mathbb{C})$, and let $\mathcal{P}$ be a totally geodesic hyperplane in $\Hh $. Let $H$ be its stabiliser in $\Gamma$. We assume that $H$ is a lattice in $\operatorname{Stab}(\mathcal{P})$,
so that $H\backslash\mathcal{P}$ has finite area. 

Define
$$ N(X;p, \mathcal{P})
:=\#\{H\gamma \in H \backslash \Gamma \mid \cosh d_{\Hh }(\gamma p,\mathcal P)\le X\}.
$$
This counts perpendicular segments of length at most
$\operatorname{arcosh}X$ from the orbit of $p$ to $\mathcal{P}$.

We write $$N(X;p, \mathcal{P})=\sum_{1<s_j \leq 2} a_j(p,\mathcal{P})X^{s_j}+E(X;p,\mathcal{P}),$$
where the coefficients $a_j(p, \mathcal{P})$ are as in equation (\ref{azp}).
Here, the sum is over the small eigenvalues $\lambda_j=s_j(2-s_j)$ of the corresponding Laplacian. 
Laaksonen \cite[Theorem 1.1]{niko} proves that $E(X;p,\mathcal{P})=O(X^{3/2})$. We provide the following improvement, which we expect to be essentially optimal.
\begin{theorem} \label{pointhyperplane}
  Let $\Gamma$, $H$, and $\mathcal{P}$ be as above. There is a set
$B_{\mathcal{P}} \subset \Gamma\backslash \Hh $ of full hyperbolic measure such that, for every $p\in B_{\mathcal{P}}$, we have
\begin{equation*}
E(X;p,\mathcal{P})
= O_{\Gamma,\mathcal{P},p}\left(
X(\log X)^{3/2}(\log\log X)^{1/2}\log\log\log{X}\right)
\end{equation*}
as $X \rightarrow \infty$.
\end{theorem}
\begin{remark} In higher dimensions, the problem generalises to counting relative to totally geodesic submanifolds. See for example \cite{Herrmann1962}. We expect our method to be applicable for that generalisation as well.
\end{remark}
We also study the spatial second moment of the error term. To that end, we further subtract a secondary term from the error term. More precisely, define $$e(X;p, \mathcal{P}):=E(X;p, \mathcal{P})-\sum_{1 < s_j < 2} b_j(p,\mathcal{P})X^{2-s_j}-c\log{X},$$
where the coefficients $b_j$ and $c$ are as in equation (\ref{bdefin}). We prove an asymptotic expansion for $$
V(X, \mathcal{P}):=\int_{\Gamma\backslash \Hh }|e(X;p,\mathcal{P})|^2 d\mu(p).$$

\begin{theorem} \label{meanvaluethm}
    For $V$ as above, as $X$ goes to infinity, we have
    $$V(X, \mathcal{P})=X^2T(X,\mathcal{P})+O_{\Gamma, \mathcal{P}}\bigl(1\bigr),$$
where, for $r\geq0$,
$$T(\cosh{r},\mathcal{P})=
4\sum_{\lambda_j\ge1}
\frac{|\widehat u_j|^2}{1+t_j^2}
\cos^2\!\bigl(t_jr-\arctan t_j\bigr).
$$
The series converges uniformly. Here, the sum in the main term is over the eigenvalues $\lambda_j=1+t_j^2 \geq 1$ of the hyperbolic Laplacian over $\Gamma \backslash \Hh$, and the $\hat{u}_j$'s are the periods of the corresponding Maa\ss{} forms over $H \backslash \mathcal{P}$, as in equation (\ref{nikoperiodsdef}).
\end{theorem}

\begin{corollary} \label{thmmvlimit} The quantity
$X^{-2}V(X, \mathcal{P})$ is bounded, but does not have a limit as $X$ goes to infinity. In particular, this implies that
    $$V(X, \mathcal{P})=O(X^2), \quad\textup{and} \quad V(X, \mathcal{P})=\Omega(X^2).$$
\end{corollary}

\begin{remark} The $O(X^2)$ result improves the upper bound $X^2\log^2{X}$ previously achieved by Laaksonen \cite[Theorem 1.3]{niko}. We should note that Laaksonen considers the initial error term $E(X;p,\mathcal{P})$, rather than $e(X;p,\mathcal{P})$. However, as $|E(X;p,\mathcal{P})-e(X;p,\mathcal{P})|=o(X)$ uniformly in $p$, this does not affect the $O(X^2)$ result.
\end{remark}

\subsection{The point-cusp problem}
In the third problem, we consider the heights of the images of a point $z \in \mathbb{H}^2$ relative to a fixed cusp $\mathfrak a$. Here, let $\Gamma$ be a cofinite, noncocompact Fuchsian group, and let $\Gamma_{\mathfrak a}$ be the stabiliser of $\mathfrak{a}$. Choose a scaling matrix $\sigma_{\mathfrak a}$ satisfying
$$
\sigma_{\mathfrak a}\infty=\mathfrak a,
\qquad
\sigma_{\mathfrak a}^{-1}\Gamma_{\mathfrak a}\sigma_{\mathfrak a}
=\langle z\mapsto z+1\rangle.
$$
We study the counting function
$$
N(X;z,\mathfrak a)
:=
\#\left\{
\gamma\in\Gamma_{\mathfrak a}\backslash\Gamma \mid
\operatorname{Im}(\sigma_{\mathfrak a}^{-1}\gamma z)\ge X^{-1}
\right\},
$$
which describes the distribution of orbit points in expanding horoballs centred at $\mathfrak a$. In the same way that the point–geodesic problem can be related to the distribution of cosets in $H\backslash\Gamma$, where $H$ is a hyperbolic cyclic subgroup, the point-cusp problem can be related to the distribution of cosets in $\Gamma_{\mathfrak{a}}\backslash\Gamma$, where $\Gamma_{\mathfrak{a}}$ is a parabolic cyclic subgroup\footnote{For an elliptic cyclic subgroup $\Gamma_1$ fixing $w$, the count
on $\Gamma_1\backslash\Gamma$ is the classical hyperbolic lattice
count centred at $w$, divided by $|\Gamma_1|$. More generally,
parabolic, hyperbolic and elliptic stabilisers give rise to counting
problems on the double-coset spaces
$\Gamma_1\backslash\Gamma/\Gamma_2$. In the point-cusp,
point-geodesic and point-point problems, $\Gamma_2$ is trivial for
generic $z$. Good \cite{Good1983} treated these problems in a
unified way, obtaining $O(X^{2/3})$ pointwise error bounds in
the appropriate counting parameters.}.

We write
$$
N(X;z,\mathfrak a)
=
M(X;z,\mathfrak{a})
+
E(X;z,\mathfrak a),
$$
where $M(X;z,\mathfrak{a})$ is the sum of the volume term $X/\hbox{vol}\left( \Gamma \backslash \H \right)$ and the contributions from exceptional poles of the parabolic Eisenstein series, as in equation (\ref{maincusp}).
The classical bound for the error term is
$$
E(X;z,\mathfrak a)=O_{\Gamma,z,\mathfrak a}(X^{2/3}).
$$
The earliest available treatment is due to Good \cite[Theorem 4]{Good1983}. A proof can be found in \cite[Theorem 4.1]{BurrinNevoRuhrWeiss2020}. We demonstrate the following almost-everywhere improvement, using the argument from \cite[Theorem 4]{Yu2018}, which is the higher-dimensional analogue (see also \cite[Theorem 4.1]{kelmer2026secondmomentincompleteeisenstein}). The proof uses a maximal estimate method due to Schmidt \cite[Theorem 2]{Schmidt1960}. For another application of this method, see \cite[Theorem 10.3]{BurrinFairchild2026}. 

\begin{theorem} \label{pointcusp} Let $\Gamma$ be a cofinite, noncocompact Fuchsian group, and fix a cusp $\mathfrak a$. There exists a set $B_{\mathfrak a}\subset\Gamma\backslash\mathbb H^2$ of full hyperbolic measure such that, for every $z\in B_{\mathfrak a}$, we have
$$
E(X;z,\mathfrak a)
=
O_{\Gamma,z,\mathfrak a}
\left(
X^{1/2}(\log X)^{3/2}
(\log\log X)^{1/2}\log\log\log X
\right)
$$
as $X\to\infty$.
\end{theorem}

 We conjecture that the exponent $1/2$ of $X$ is optimal for the generic case. To that end, we prove the following result, using an argument similar to \cite[Theorem 3.3]{philipsr}.
 \begin{theorem} \label{omegatheorem}
     If $\Gamma$ is not essentially cuspidal, then there exists a cusp $\mathfrak{a}$ such that, for almost every $z$, and every $\delta>0$, we have $E(X;z,\mathfrak{a})=\Omega(X^{1/2-\delta})$.
 \end{theorem}
 Here, a cofinite Fuchsian group $\Gamma$ is called \emph{essentially cuspidal} if $$\#\left\{j \mid 0\leq \lambda_j \leq T\right\} \sim \frac{\operatorname{vol}(\Gamma\backslash\mathbb H^2)}{4 \pi}T,$$ where the sequence $(\lambda_j)_j$ denotes the eigenvalues of the discrete spectrum of the corresponding Laplacian (see \cite[equation (25)]{Sarnak2003}). We note that, while the existence of groups that are not essentially cuspidal is not known, it is widely believed that, in a certain sense, most groups should not be essentially cuspidal (see \cite[Conjecture 1]{Sarnak2003}, \cite{PhillipsSarnak1985}, \cite{Wolpert1994}, and \cite[p.160]{Iwaniec2002}). Furthermore, their existence is also known conditional on an eigenvalue multiplicity hypothesis (see \cite[Theorem 2]{luo2001}).
 
It is worth mentioning that the conclusion of Theorem 6 does not hold for every cofinite group. In particular, for the modular group, the count can be expressed as a count of primitive lattice points in an ellipse, and the quasi-Riemann hypothesis \cite{OpenAI2026QRH, ji2026proofquasiriemannhypothesis} yields a power saving over $X^{1/2}$ by a standard M\"obius-inversion argument similarly with the proof of \cite[Theorem 13]{Risager2026}. See \cite[\S 1, pp. 520-521]{Nowak2005} for the connection with zero-free regions of $\zeta(s)$. Conditional to the full Riemann hypothesis, the best known upper bound for the error term is $X^{5/13+\epsilon}$ (see \cite{RogerC2010}).

\subsection{Structure} In Section \ref{sectionmaxest}, we formulate the maximal estimate method for discrete and continuous indexing. Section \ref{secth2} contains the proofs of the point–geodesic and point–cusp results in $\mathbb H^2$, and Section \ref{secth3} treats the point–hyperplane problem in $\mathbb H^3$, including the spatial second-moment formula.
\section{Maximal estimates} \label{sectionmaxest}
Let $(\mathcal X,\mu)$ be a measure space and let
$S\colon [1,\infty)\times \mathcal{X}\to\mathbb C$.
By a maximal estimate for $S$ on an interval
$I\subseteq [1,\infty)$, we mean a bound for
$$
\int_{\mathcal{X}} \sup_{R\in I}|S(R,z)|^2\,d\mu(z).
$$
Applied on successive intervals, such estimates can be combined with
Chebyshev's inequality and the Borel--Cantelli lemma to obtain pointwise bounds for $S(R,z)$
that hold for almost every $z$ and all sufficiently large $R$.
This is the key idea behind what we refer to as a maximal estimate method. The central question is therefore how to establish these maximal estimates.

In \cite{chatzakos2026squarerootcancellationhyperbolic}, Chatzakos--Dimakis use such a method to prove square-root cancellation in the hyperbolic lattice counting problem, for almost every base point $z$.
Essentially, their method relies on two observations. First, grouping the spectral terms into dyadic blocks lets us exploit averaged spectral estimates, such as the local Weyl law, at each scale. Second, bounds for the $R$-derivative of a block give information about its variation over short intervals. A Sobolev inequality combines this information with a bound for the block itself to control its supremum on each interval. We expect the method to be useful in other counting problems where comparable mean-square bounds for the blocks and their derivatives are available.

We note that other maximal estimate methods exist as well. In \S \ref{pointcuspsec}, we use a different maximal argument for the point–cusp problem, based on mean-square bounds for increments and a dyadic decomposition in the counting parameter. For an earlier
use of this dyadic argument, see \cite[Lemmas~2 and~3]{Schmidt1960}.

In this section, we formulate the Chatzakos--Dimakis argument in a general setting,
first for discrete indexing and then for continuous indexing.
\subsection{Discrete indexing}
Let $(\mathcal{X}, \mu)$ be a measure space. Let $(\phi_j)_j$ be a sequence of orthonormal $L^2-$functions defined on $\mathcal{X}$, and let $c_j(r)$ be a sequence of functions in $C^1([1,\infty))$.

We want to find an upper bound for \begin{equation*} \label{srzdefinition}S(R,z):=\sum_{j}c_j(R)\phi_j(z),\end{equation*}
in terms of $R$. We assume that, for every $z \in \mathcal{X}$ and every compact interval $I\subseteq[1,\infty)$, the series defining $S(R,z)$ converges uniformly for $R \in I$. We call this condition (C0). This can be established, for example, by a local Weyl law on the $\phi_j$'s, paired with an appropriate decay condition for the coefficients.

Suppose that there exists a sequence $(t_j)_j$ of real numbers with $$1 \leq t_1 \leq t_2 \leq \dots \leq t_j \leq \dots$$ and $t_j \rightarrow \infty$, and a function $D(t,r)$ taking positive real values only, with the following properties:
\begin{itemize}
\item[(C1)] For any real numbers $T,r \geq 1$, we have $$ \sum_{T \leq t_j <2T}\left(|t_jc_j(r)|^2+|\partial_rc_j(r)|^2\right) \ll  T(D(T,r))^2.$$
    \item[(C2)] For any real numbers $t,a \geq 1$, we have $$ \sup\limits_{r \in [a,a+1]}D(t,r) \ll D(t,a).$$
    \item[(C3)]  The series $s(x):=\sum\limits_{k=0}^{\infty}D(2^k,x)$ converges for every $x \geq 1$.
\end{itemize}

Then, we have the following theorem:
\begin{theorem} \label{mainmaximal}
     Let $g:[1, \infty) \rightarrow (0,\infty)$ be a nondecreasing function such that $$\sum \left(\frac{s(n)}{g(n)}\right)^2 < \infty.$$
    Then, for $\mu$-almost every $z\in  \mathcal{X}$, we have
$$
|S(R,z)|\ll_{z,g} g(R).
$$
\end{theorem}
\begin{remark}
If the partial sums of $S(R,z)$ over $j$ do not converge, we can define instead $S(R,z)$ by
$$S(R,z):=\lim_{K \rightarrow \infty} S_K(R,z),$$ when this converges, and $0$ otherwise.
Here, $$S_K(R,z):=\sum_{t_j<2^{K+1}}c_j(R)\phi_j(z).$$
Using Proposition \ref{keystep} and (C3), one can prove that this is guaranteed to converge locally uniformly in $R$ for $z$ lying in a full-measure set independent of $R$. In that case, the theorem still holds when further restricted to the set of convergence.
\end{remark}
First, we introduce a dyadic decomposition in the spectral parameter. Since there are finitely many $t_j$'s in $[T,2T]$, this will give us finite sums of functions, which are easier to work with than infinite sums. 

We define 
$$B(T,R,z):=\sum_{T\leq t_j<2T}c_j(R)\phi_j(z).$$
Using the orthogonality of the functions $\phi_j$, we have the following $L^2$-estimate with respect to $z$.
 \begin{lemma} \label{bestimate}

We have
$$T\|B(T,R, \cdot)\|_2^2,\, \,\,T^{-1}\|\partial_R B(T,R,\cdot )\|_2^2 \ll D^2(T,R).$$
 \end{lemma}
 \begin{proof}
     Using orthogonality and condition (C1), we have $$T^2\|B(T,R,\cdot)\|_2^2= T^2\sum_{T\leq t_j<2T}|c_j(R)|^2 \ll \sum_{T\leq t_j<2T}|t_jc_j(R)|^2 \ll TD(T,R)^2.$$
     We proceed similarly for $\partial_RB(T,R,\cdot )$.
 \end{proof}

 The next step is to bound the $L^2$ $z$-norm of $\sup_{R \in I}|B(T,R,z)|$, where $I$ is an interval of unit length.  A Sobolev inequality in $\R$, followed by Cauchy–Schwarz and Tonelli, reduces this to the two estimates in Lemma \ref{bestimate}.
 This step is the key idea of the proof.

 \begin{proposition} \label{keystep} For all real numbers $a,T \geq 1$, we have
     $$\int \sup\limits_{R \in [a,a+1]}|B(T,R,z)|^2 d\mu(z) \ll D^2(T,a).$$
 \end{proposition}
 \begin{proof}
 Let $I=[a,a+1]$. We write $\| \, \cdot \,\|_{k,I}$ for the $L^k$ norm over $I$.
 
Using the Sobolev inequality of Chatzakos--Dimakis \cite[Lemma 3.4]{chatzakos2026squarerootcancellationhyperbolic}, we have
 $$\|B(T,\, \cdot \,,z)\|_{\infty,I}^2 \ll \|B(T,\, \cdot \,,z) \|_{2,I}^2+\|B(T,\, \cdot \,,z) \|_{2,I}\|\partial_RB(T,\, \cdot \,,z) \|_{2,I}. $$
Hence,

 $$\int \|B(T,\, \cdot \,,z)\|_{\infty,I}^2 \, \, d\mu(z) \ll \int \|B(T,\, \cdot \, ,z)  \|_{2,I}^2d\mu(z)+ \int \|B(T,\, \cdot \, ,z) \|_{2,I}\|\partial_RB(T,\, \cdot \, ,z) \|_{2,I} d\mu(z).$$
 Using Cauchy--Schwarz, this gives
 \begin{align*} & \int \sup\limits_{R \in I}|B(T,R,z)|^2 d\mu(z) \\ 
 &\ll\int \|B(T,\, \cdot \,,z)  \|_{2,I}^2d\mu(z)+ \left( \int \|B(T,\, \cdot \,,z) \|_{2,I}^2d\mu(z) \right)^{1/2} \left( \int \|\partial_RB(T,\, \cdot \,,z) \|^2_{2,I} d\mu(z)\right)^{1/2}. \end{align*}

 By Tonelli’s theorem, we can exchange the order of integration in the right hand side, which gives
 \begin{align*} & \int \sup\limits_{R \in I}|B(T,R,z)|^2 d\mu(z) \\ 
 &\ll \int_{I} \|B(T,R, \cdot)\|_2^2\,dR + \left(\int_{I} \|B(T,R, \cdot)\|^2_2\,dR \right)^{\frac{1}{2}}\left(\int_{I} \|\partial_RB(T,R, \cdot)\|_2^2 dR\right)^{\frac{1}{2}}. \end{align*}

 Hence, applying Lemma \ref{bestimate} and condition (C2), we have
$$\int \sup\limits_{R \in I}|B(T,R,z)|^2 d\mu(z) \ll T^{-1}|D(T,a)|^2+ |D(T,a)|^2 \ll |D(T,a)|^2. $$
 
 \end{proof}

  Note that $$S(R,z):=\sum_{k=0}^{+\infty}B(2^k,R,z).$$
We will now use the Minkowski inequality to combine the dyadic estimates  in $S(R,z)$. For every $a \geq 1$, we define
 $$G_a(z)=\sup_{R \in [a,a+1]}|S(R,z)|.$$
 We have the following upper bound for its second moment.

 \begin{lemma} \label{glemma} For every $a \geq 1$, we have
    $$ \int |G_a(z)|^2 d \mu (z) \ll (s(a))^2,$$
    where $$s(x):=\sum\limits_{k=0}^{\infty} D(2^k,x).$$
 \end{lemma}
 \begin{proof}
     We write $S_K(R,z)$ for the truncated sum up to $k=K$. Similarly, we write $$G_{a,K}(z)=\sup_{R \in [a,a+1]}|S_K(R,z)|.$$
  Applying the Minkowski inequality and Proposition \ref{keystep}, we obtain 
\begin{equation*}
\begin{split}
\left(\int |G_{a,K}(z)|^2\,d\mu(z)\right)^{\frac{1}{2}} &= \left(\int \sup\limits_{R\in [a,a+1]} \left|\sum \limits_{k=0}^{K} B(2^k,R,z)\right|^2\,d\mu(z) \right)^{\frac{1}{2}}\\
& \leq  \sum\limits_{k=0}^K \left(\int\sup\limits_{R\in [a,a+1]}\left|B(2^k,R,z)\right|^2\,d\mu(z)\right)^{\frac{1}{2}}\\
&\ll \sum\limits_{k=0}^K D(2^k,a) \\
&\ll s(a).
\end{split}
\end{equation*}
 Finally, by uniform convergence and Fatou's lemma, we have
\begin{equation*}
\int|G_a(z)|^2\,d\mu(z) \le \liminf\limits_{K\to \infty} \int|G_{a,K}(z)|^2\,d\mu(z) \ll (s(a))^2. 
\end{equation*}
 \end{proof}
Fix a nondecreasing positive real function $g(x)$. Consider the set of points  
\begin{equation*}
E_{n,g} = \{z \mid G_n(z) > g(n)\}. 
\end{equation*}
The Chebyshev inequality and Lemma \ref{glemma} imply that 
\begin{equation*}
\mu(E_{n,g})  \leq \frac{1}{g^2(n)} \int |G_n(z)|^2\,d\mu(z) \ll \left(\frac{s(n)}{g(n)}\right)^2. 
\end{equation*}
Assume now that the series
$$\sum \left(\frac{s(n)}{g(n)}\right)^2$$ converges. Then,
\begin{equation*}
\sum\limits_{n=1}^{\infty} \mu(E_{n,g}) <\infty. 
\end{equation*}
Thus, the first Borel-Cantelli Lemma implies that 
\begin{equation*}
\mu(\limsup\limits_{n\to \infty} (E_{n,g})) = 0. 
\end{equation*}
Define $B_g:=  \mathcal{X}\backslash \limsup\limits_{n\to \infty} (E_{n,g})$. If $z\in B_g$ then there exists an $n_0= n_0(z,g)$ such that 
\begin{equation*}
G_n(z) \leq g(n) \text{ for } n\geq n_0. 
\end{equation*}
For $R\in [n,n+1]$ this implies that 
\begin{equation*}
|S(R,z)|\leq g(n) \leq g(R).
\end{equation*}
This concludes the proof of Theorem \ref{mainmaximal}.
\subsection{Continuous indexing}

In this subsection we present a version of Theorem \ref{mainmaximal} with continuous indexing.
As before, let $(\mathcal{X},\mu)$ be a measure space.  

Let $\Phi:\mathcal{X} \times \mathbb R\to\mathbb C$ be jointly
measurable, and suppose that, for every $z\in  \mathcal{X}$,
the function $t\mapsto\Phi(z,t)$ is continuous on $\mathbb R$.
The function $\Phi(z,t)$ is the analogue of the sequence $\phi_j(z)$.  

Let $c:[1,\infty)\times\mathbb R\to\mathbb C$ be continuous,
and assume that $\partial_Rc$ exists and is jointly continuous. The function $c(R,t)$ is the analogue of the sequence $c_j(R)$.

Furthermore, we assume that,
uniformly for $T,R \geq 1$, we have
$$
\left\|
\int_{ T}^{2T} a(R, t)\Phi(\cdot, t)\,dt
\right\|_{L^2( \mathcal{X})}^2
\ll\int_{T}^{2T} |a(R, t)|^2\,dt,
$$
for both $a(R,t)=c(R,t)$ and $a(R,t)=\partial_Rc(R,t)$. This Bessel bound substitutes the orthonormality property we had for the sequence $\phi_j$.

 Define
$$
S_K(R,z):=
\int_{1}^{2^K}c(R,t)\Phi(z,t)\,dt.
$$

We assume that, for every $z\in \mathcal{X}$, the sequence $S_K(R,z)$ converges locally uniformly for $R\geq1$. We call this assumption (C0)', as it is the analogue of (C0). Denote its limit by $S(R,z)$.

Furthermore, we assume that there is a positive function $D(T,R)$ that satisfies properties (C2),(C3), and the following analogue of property (C1):
\begin{itemize}
\item[\textup{(C1)'}] Uniformly for $T,R\ge1$,
  $$
  \int_{T}^{2T}
  \left(t^2|c(R,t)|^2+
  |\partial_Rc(R,t)|^2\right)
  dt
  \ll T D(T,R)^2.
  $$
  \end{itemize}

Then, we have the following analogue of Theorem \ref{mainmaximal}:
\begin{theorem} \label{contmainmaximal}

Let $g:[1,\infty)\to(0,\infty)$ be nondecreasing and satisfy
$$
\sum_{n=1}^{\infty}\left(\frac{s(n)}{g(n)}\right)^2<\infty.
$$
Then, for $\mu$-almost every $z\in  \mathcal{X}$, we have
$$
|S(R,z)|\ll_{z,g} g(R).
$$
\end{theorem}

\begin{proof}
    Let $$B(T,R,z):=\int_{ T}^{2T} c(R, t)\Phi(z,t)\,dt.$$
    Using property (C1)', we have $$T\|B(T,R, \cdot)\|_2^2,\, \,\,T^{-1}\|\partial_R B(T,R,\cdot )\|_2^2 \ll D^2(T,R),$$
    as in Lemma \ref{bestimate}. The proof then follows as in the previous section, essentially unchanged.

\end{proof}

\section{Problems in $\H$} \label{secth2}
Let $\Gamma$ be a cofinite Fuchsian group, and let $(u_j)_j$ be a maximal orthonormal set of real-valued eigenfunctions for the discrete spectrum of the Laplacian in $L^2(\Gamma\backslash\mathbb H^2,d\mu)$, with corresponding eigenvalues $\lambda_j=s_j(1-s_j)$. We write $s_j=1/2+it_j$, where $t_j$ is either nonnegative real or purely imaginary. More precisely, either  $\hbox{Re}(s_j)=1/2$ and $\hbox{Im}(s_j) \geq 0$, or $s_j \in [1/2,1]$.

For every cusp $\mathfrak{b}$, let $\Gamma_{\mathfrak b}$ be its stabiliser. Choose a scaling matrix $\sigma_{\mathfrak b}$ satisfying
$$
\sigma_{\mathfrak b}\infty=\mathfrak b,
\qquad
\sigma_{\mathfrak b}^{-1}\Gamma_{\mathfrak b}\sigma_{\mathfrak b}
=\langle z\mapsto z+1\rangle.
$$
Let
$$
\mathcal{E}_{\mathfrak b}(z,s)
=
\sum_{\gamma\in\Gamma_{\mathfrak b}\backslash\Gamma}
\operatorname{Im}(\sigma_{\mathfrak b}^{-1}\gamma z)^s,
\qquad \operatorname{Re}s>1
$$
be the associated Eisenstein series, and let $\mathcal P_{\mathfrak b}\subset(1/2,1)$ denote its exceptional poles.

While the growth of individual eigenfunctions $u_j$ and $\mathcal{E}_{\mathfrak{b}}(z,s)$ is hard to study in the spectral aspect, we have the following mean square result.
\begin{lemma}[Local Weyl law {\cite[Proposition 7.2]{Iwaniec2002}} ] \label{lemmaweyl} For $T \geq 1$, we have
$$
\sum_{|t_j|\le T}|u_j(z)|^2
+
\sum_{\mathfrak b}
\int_{-T}^{T}
\left|\mathcal{E}_{\mathfrak b}\left(z,\tfrac{1}{2}+it\right)\right|^2\,dt
\ll_{\Gamma,z}T^2.
$$
\end{lemma}
\subsection{The point-geodesic problem}
\subsubsection{Background theory}
Let $l= H \backslash \tilde{l}$ be a fixed oriented closed geodesic with stabiliser $H$. By conjugating the group, we assume that $\tilde{l}$ is the positive imaginary axis, which we denote by $I$. We denote by $\hbox{len}(l)$ the hyperbolic length of $l$:
$$\hbox{len}(l):=\int_{l} 1 \; ds(z).$$
For $z=x+iy \in \H$, we define the \emph{Huber coordinates} $(u,v)$ as follows:
\begin{equation*} \label{coordsys} u(z)=\log{|z|}, \; \; \; \; v(z)=-\arctan{\left(\frac{x}{y}\right)},\end{equation*}
or, equivalently,
$$x=-e^u\sin{v}, \; \; \; \; y=e^u\cos{v}. $$  
We note that, for $\gamma$ diagonal, we have $$u\left(\gamma z \right)=u(z)+\log{\lambda},$$ where $\lambda$ is the norm of $\gamma$, and
$$v\left(\gamma z\right)=v(z).$$
The quantity $v(z)$ can be interpreted as the anticlockwise angle formed between $z=ie^{u+iv}$ and the positive imaginary axis.
With respect to Huber coordinates, 
\begin{equation*}\label{meshub} \, d\mu(z)=\frac{1}{\cos^2{v}}\, du\, dv. \end{equation*}
Furthermore, we define the \emph{periods} of eigenfunctions and Eisenstein series by
\begin{equation*} \label{periodsdef}\hat{u}_{j}:=\int_{l}u_{j}(z) ds(z), \quad \hat{\mathcal{E}}_{\mathfrak{b}}(s'):=\int_{l}\mathcal{E}_{\mathfrak{b}}\left(z, \, s'\right) ds(z).\end{equation*}

For a bounded, compactly supported $f$ on $[1,\infty)$, Huber considers the following series:
\begin{equation} \label{hubers} A(f)(z):=\sum_{\gamma \in H \backslash \Gamma}f\left(\frac{1}{\cos^2{\left( v(\gamma z)\right)}}\right). \end{equation}
This is an automorphization of $f$ composed with $1/\cos^2{\left( v(z)\right)}$, with respect to $H \backslash \Gamma$, where $1/\cos^2{\left( v(z)\right)}$ is associated to the hyperbolic distance of $z$ from the imaginary axis, $\rho(z,I)$, via the relation $1/\cos^2{\left( v(z)\right)}=\cosh^2{\rho(z,I)}$. Therefore, for $f(x)$ the indicator function of $\left\{ x \leq X^2 \right\}$,
we have that
$A(f)(z)=N(X;z,l). $
Huber \cite{Huber1998}, and later Chatzakos--Petridis \cite{ChatzakosPetridis+2016+981+1003}, take $f$ to be a smooth approximation of the indicator, and consider the spectral expansion of the function $A(f)(z)$. In particular, they prove the following result. \begin{lemma}[{\cite[Lemmas 2.1 and 4.2]{ChatzakosPetridis+2016+981+1003}}] \label{a0spectral}
For $f$ a continuous, compactly supported, piecewise $C^2$ function,
we have the following spectral expansion:
$$A(f)(z)=2\sum_{j}d(f,t_j)\hat{u}_{j}u_{j}(z)-\sum_{\mathfrak{b}}\frac{i}{2\pi}\int_{\left(1/2\right)}d(f,t)\hat{\mathcal{E}}_{\mathfrak{b}}\left(s\right)\overline{\mathcal{E}_{\mathfrak{b}}\left(z,s\right)}ds,$$
where
\begin{equation*}
d(f,t)=\int_0^{\pi/2}
f(\cos^{-2}v)
\frac{\xi_{1/4+t^2}(v)}{\cos^2 v}\,dv,
\end{equation*}
and
\begin{equation} \label{difequation}
\xi_\lambda''(v)+\lambda\cos^{-2}v\,\xi_\lambda(v)=0,
\qquad
\xi_\lambda(0)=1,
\qquad
\xi_\lambda'(0)=0.
\end{equation}
Here, piecewise $C^2$ assumes a finite partition on the support, with finite one-sided limits of the derivatives through order two.
\end{lemma}

To deal with the contribution of the periods, they use the following estimate of Huber.

\begin{lemma}[{Huber \cite{Huber1998}, Chatzakos--Petridis \cite[Lemma 4.3]{ChatzakosPetridis+2016+981+1003}}]\label{lemmahuber} For $T \geq 1$, we have
$$\sum_{|t_j|\leq T} \left| \hat{u}_{j} \right|^2
+ \sum_{\mathfrak{b}}  \int_{-T}^{T} \left|\hat{\mathcal{E}}_{\mathfrak{b}}\left(\tfrac{1}{2}+it \right) \right|^2 dt 
\ll_{\Gamma, l} T
    . $$

\end{lemma}

By estimating carefully the Huber transforms $d(f,t)$ that arise, and using Lemma \ref{lemmaweyl} and Lemma \ref{lemmahuber} for the values of $u_j$ and the periods respectively, they show that, for fixed $z \in \mathbb{H}^2$, we have that 
    $$N(X;z,l)=\sum_{1/2 < s_j \leq 1} a_j(z,l)   X^{s_j}+O(X^{2/3}),$$
where
\begin{equation} \label{ajzldef} a_j(z,l)=A(s_j)\hat{u}_{j} u_{j}(z) \end{equation}
and
$$
A(s)=
\frac{2\sqrt\pi\,\Gamma(s-\tfrac12)}
{s\,\Gamma(s/2)^2}.
$$
\subsubsection{Proof of Theorem \ref{pointline}}

\label{sectionproofpointline}

Put $U=(X^2-1)^{1/2}$, and let $Y$ be a real number with $0<Y<U/2$.
We define the \emph{upper comparison function} $f^+_{U,Y}$ by
\begin{equation*}
f_{U,Y}^+(x^2+1)=
\begin{cases}
1, & 0\leq x\leq U,\\
(U+Y-x)/Y, & U\leq x\leq U+Y,\\
0, & x\geq U+Y,
\end{cases}
\end{equation*}
and the \emph{lower comparison function} $f^-_{U,Y}$ by
\begin{equation*} 
f_{U,Y}^-(x^2+1)=
\begin{cases}
1, & 0\leq x\leq U-Y,\\
(U-x)/Y, & U-Y\leq x\leq U,\\
0, & x\geq U.
\end{cases}
\end{equation*}
Note that $f_{U,Y}^{-}=f^+_{U-Y,Y}$. Furthermore, we write 
$$f^{\pm}_{U}(x):=f^{\pm}_{U,1}(x).$$
Then,
\begin{equation}  \label{compareerror}
A(f_{U}^-)(z)\leq N(X;z,l)\leq A(f_{U}^+)(z).
\end{equation}
By Lemma \ref{a0spectral}, we have

$$A(f_{U}^{\pm})(z)=2\sum_{j}d(f_U^{\pm},t_j)\hat{u}_{j}u_{j}(z)-\sum_{\mathfrak{b}}\frac{i}{2\pi}\int_{\left(1/2\right)}d(f_U^{\pm},t)\hat{\mathcal{E}}_{\mathfrak{b}}\left(s\right)\overline{\mathcal{E}_{\mathfrak{b}}\left(z,s\right)}ds.$$
Defining $$E^{\pm}(X,z):=A(f^{\pm}_U)(z)-\sum_{1/2 < s_j \leq 1} a_j(z,l)   X^{s_j},$$
we have that it is enough to prove the required bound for $E^{\pm}$. We will do this by applying
Theorems \ref{mainmaximal} and \ref{contmainmaximal}. To do this, we first provide some estimates for the corresponding transforms
and their derivatives.

\begin{lemma} \label{dfutest}
Write $X=e^r$, $U=(X^2-1)^{1/2}$. For $r \geq1$, we have the following estimates:
\begin{itemize}
\item[\textup{(a)}]
For $0 \leq t \leq 1$, we have $$|d(f_{U}^{\pm},t)| \ll X^{1/2}\log{X}.$$
\item[\textup{(b)}]
For fixed $s=1/2+it \in (1/2,1]$, in which case $t$ is purely imaginary, we have $$2d(f_{U}^{\pm},t)=A(s)X^s+O_s(X^{1/2}).$$
    \item[\textup{(c)}]
For $t\geq 1$, set
\begin{equation} \label{defqtr1}
Q(t,r)=\min\{1,X/t\}.
\end{equation}

Then, uniformly in $r$ and $t$,
\begin{equation} \label{dfestimate}
X^{-1/2}|d(f_{U}^{\pm},t)|
\ll t^{-1}Q(t,r),
\end{equation}
and
\begin{equation}
\left|
\partial_r\left(
X^{-1/2}d(f_{U}^{\pm},t)
\right)
\right|
\ll Q(t,r).
\end{equation}

\end{itemize}
\end{lemma}

\begin{proof}
We demonstrate this for $f_{U}^+$. The case of $f_{U}^-$ is similar.
For part (a), we note that, from \cite[p.9 \& equation (6)]{lekkas2024}, it follows that there is a nonnegative kernel $k$ such that
$d(f_{U}^{+},t)=h_k(t)$, where $h_k$ is the Selberg--Harish--Chandra transform of $k$.
Hence,
$$
|d(f_{U}^{+},t)|
=\left|\int_{\mathbb H} k(u(i,z))\,y^{1/2+it}\,d\mu(z)\right|
\le \int_{\mathbb H} k(u(i,z))\,y^{1/2}\,d\mu(z)
=d(f_{U}^{+},0).
$$

Therefore, from \cite[Proposition 2.4(b),(d)]{ChatzakosPetridis+2016+981+1003}, we have $|d(f_{U}^{+},t)|\ll X^{1/2}\log{X}$. This proves part (a). Part (b) follows directly from \cite[Proposition 2.4(c)]{ChatzakosPetridis+2016+981+1003}.

Furthermore, from \cite[Proposition 2.4(b)]{ChatzakosPetridis+2016+981+1003},  we
have 
\begin{equation} \label{dfesty}
X^{-1/2}|d(f_{U,Y}^{+},t)|
\ll t^{-1}\min\{1,X/(tY)\}.
\end{equation}
The estimate (\ref{dfestimate}) follows directly by setting $Y=1$.
Now write
$$
D_t(x)=\int_0^x\xi_{1/4+t^2}(\arctan s)\,ds.
$$
Taking $Y \rightarrow 0^+$ in equation (\ref{dfesty}), we get 
$$
|D_t(x)|=|\lim_{Y \rightarrow0}d(f^+_{x,Y},t)|\ll t^{-1}(1+x^2)^{1/4}.
$$

On the other hand, writing
$\lambda=t^2+\tfrac14$, equation (\ref{difequation}) becomes
$$
\bigl((1+x^2)D_t''(x)\bigr)'+\lambda D_t'(x)=0.
$$
Since $D_t(0)=D_t''(0)=0$, integration gives
$$
D_t''(x)=-\frac{\lambda}{1+x^2}D_t(x).
$$
Choose $h=\sqrt{1+x^2}/t$. Since $t\ge1$, we have $1+s^2\asymp1+x^2$ throughout $s \in [x,x+h]$. Taylor’s formula therefore yields
$$
\begin{aligned}
|D_t'(x)|
&\le
\frac{|D_t(x+h)|+|D_t(x)|}{h}
+\frac h2\sup_{x\le s\le x+h}|D_t''(s)|\\
&\ll
\frac{(1+x^2)^{1/4}}{th}
+ht(1+x^2)^{-3/4}\\
&\ll(1+x^2)^{-1/4}.
\end{aligned}
$$
Hence,
\begin{equation} \label{Destimates}|tD_t(x)| \ll (1+x^2)^{1/4}, \; \; \hbox{and} \quad \, |D_t'(x)| \ll (1+x^2)^{-1/4}. \end{equation}

By Fubini, we have
$$d(f^+_{U},t)=\int_{U}^{U+1}D_t(s) ds.$$

Differentiating this identity gives
$$
\partial_r d(f^+_{U},t)
=U'\bigl(D_t(U+1)-D_t(U)\bigr)
$$
On the other hand, using equation (\ref{Destimates}) and the mean value theorem, we have
$$
|D_t(U+1)-D_t(U)|
\ll \min\{X^{-1/2},X^{1/2}/t\}.
$$
Hence,
\begin{align*}
|\partial_r d(f^+_{U},t)|
\ll U' \cdot\min\{X^{-1/2},X^{1/2}/t\}
\ll X \cdot\min\{X^{-1/2},X^{1/2}/t\}
= X^{1/2} Q(t,r).
\end{align*}
This gives

$$|\partial_r (X^{-1/2}d(f^+_{U},t)) |\ll X^{-1/2} \left(|d(f^+_{U},t)|+ |\partial_r d(f^+_{U},t)|\right)
\ll  Q(t,r).$$
\end{proof}
We now define 
$$E_{\text{disc.}}^{\pm}(X,z):= 2\sum_{t_j \geq 1} d(f^{\pm}_{U},t_j)\hat{u}_ju_j(z)$$
and, for every class of cusps $\mathfrak{b}$,

$$E_{\mathfrak{b}}^{\pm}(X,z):=
\hbox{Re}\int_{1}^{\infty}
d(f_U^\pm,t)\,
\overline{\widehat{\mathcal{E}}_{\mathfrak{b}}(\tfrac12+it)}
\,\mathcal{E}_{\mathfrak{b}}(z,\tfrac12+it)\,
\frac{dt}{\pi}. $$
Using Lemma \ref{dfutest}(a),(b), we have

$$E^{\pm}(X,z)=E_{\text{disc.}}^{\pm}(X,z)+\sum_{\mathfrak{b}}E^{\pm}_{\mathfrak{b}}(X,z)+O(X^{1/2}\log{X}).$$

Therefore, it is enough to prove the required bound for $E_{\text{disc}}^{\pm}$ and $E^{\pm}_{\mathfrak{b}}$.

 We start with $E_{\text{disc}}^{\pm}$. We will use Theorem \ref{mainmaximal}. If the discrete spectrum is finite, then Lemma \ref{dfutest}(c) implies $E_{\text{disc}}^{\pm}=O_{\Gamma,z,l}(X^{1/2})$, for every $z$. Suppose now that the discrete spectrum is infinite. We let $t_j$ be the spectral parameters that are $\geq 1$, and $\phi_j=u_j$. Furthermore, we let $$c_j(r)=2e^{-r/2}d(f_{U}^{\pm},t_j)\hat{u}_j, \; \; \hbox{and} \quad D(T,r)=Q(T,r).$$

 We will show that this choice satisfies the conditions (C0), (C1), (C2), and (C3).
 \begin{lemma} \label{conditioncheck}
     For $t_j$, $c_j$, $\phi_j$ and $D(t,r)$ as above, the conditions \textup{(C0), (C1), (C2), (C3)} are satisfied. Moreover, for $r$ large enough,
     \begin{equation} \label{sx1} s(r):=\sum_{k \geq 0}D(2^k,r)\ll r. \end{equation}
 \end{lemma}
 \begin{proof}

For condition (C0), let $I \subset [1,\infty)$ be a compact interval. We have

$$
\begin{aligned}
\sum_{T\le t_j<2T}\sup_{r\in I}|c_j(r)u_j(z)|
&\ll_I T^{-2}\sum_{T\le t_j<2T}|u_j(z)\hat{u}_j| \\
&\ll T^{-2}
\left(\sum_{T\le t_j<2T}|\widehat u_j|^2\right)^{1/2}
\left(\sum_{T\le t_j<2T}|u_j(z)|^2\right)^{1/2}\\
&\ll T^{-1/2},
\end{aligned}
$$
where we used Lemma \ref{lemmahuber} and Lemma \ref{lemmaweyl}. Summing over $T=2^k$ proves absolute local uniform convergence, and hence (C0).
 
     For condition (C1), we have
    \begin{align*}\sum_{T \leq t_j <2T}\left(|t_jc_j(r)|^2+|\partial_rc_j(r)|^2\right) &\ll \sum_{T \leq t_j <2T}\left(|t_j^2t_j^{-2}|D(t_j,r)|^2|\hat{u}_j|^2+|D(t_j,r)|^2|\hat{u}_j|^2\right) \\
    &\ll     |D(T,r)|^2 \sum_{T \leq t_j < 2T} |\hat{u}_j|^2. \end{align*}
    Hence, using Lemma \ref{lemmahuber}, we have
    $$\sum_{T \leq t_j <2T}\left(|t_jc_j(r)|^2+|\partial_rc_j(r)|^2\right) \ll     T|D(T,r)|^2,$$
    as required.

    Condition (C2) follows directly from equation (\ref{defqtr1}).

    For condition (C3), we have 
    \begin{equation*} \sum_{k=0}^{K}D(2^k,r)= \sum_{2^k  \leq \min\{e^r,2^K\}}1+\sum_{e^r < 2^k \leq 2^K}\frac{e^r}{2^k}\ll r+ 1 \ll r, \end{equation*}
    which is independent of $K$.
    Therefore, the series defining $s(r)$ converges for every $r \geq 1$.
 \end{proof}

From the above, it follows that we can apply Theorem \ref{mainmaximal}, for any admissible function $g$. 
We note that, by equation (\ref{sx1}), we have
$s(n) \ll n$. Hence, setting $g(x) =x^{3/2}(\log{x})^{1/2}\log{\log{x}}$ for $x \geq 4$ and $g(x)=g(4)$ otherwise, we have
\begin{equation} \label{gadmissible} \sum_{n} \left( \frac{s(n)}{g(n)}\right)^2 \ll\sum_{n \geq 4} \frac{1}{n \log{n} (\log\log{n})^2 } < \infty. \end{equation}
Therefore, there are sets $B_{\pm} \subset \Gamma \backslash \H$ of full measure, such that, for every $z \in B_{\pm}$, we have $$E_{\text{disc.}}^{\pm}(X,z)=O_{\Gamma,z,l}(X^{1/2}g(\log{X})).$$

For the continuous spectrum, fix a class of cusps $\mathfrak{b}$. For $t \geq 1$, we take $$c(r,t):=
\frac1\pi e^{-r/2}
d(f_U^\pm,t)\,
\overline{\widehat{ \mathcal{E}}_{\mathfrak{b}}(\tfrac12+it)}, \quad D(t,r)=Q(t,r) \quad \hbox{and} \quad \Phi(z,t)=
\mathcal{E}_{\mathfrak{b}}(z,\tfrac12+it).$$ 
We have already established conditions (C2),(C3). Conditions (C0)' and (C1)' can be established in the same way as with conditions (C0) and (C1). Finally, the Bessel bound condition follows from \cite[Proposition 7.1]{Iwaniec2002}.

Applying Theorem \ref{contmainmaximal} for each individual class of cusps $\mathfrak{b}$, we have that there are sets $B^{\mathfrak{b}}_{\pm}$ of full measure, such that, for every $z \in B^{\mathfrak{b}}_{\pm}$, we have $$E_{\mathfrak{b}}^{\pm}(X,z)=O_{\mathfrak{b},\Gamma,z,l}(X^{1/2}g(\log{X})).$$
Let $$B=(B_{+} \cap B_{-})\bigcap_{\mathfrak{b}}(B_{+}^{\mathfrak{b}}\cap B_{-}^{\mathfrak{b}}).$$
From equation (\ref{compareerror}), it follows that, for every $z \in B$, we have
$$E(X;z,l)=O_{\Gamma,z,l}(X^{1/2}g(\log{X}))=O_{\Gamma,z,l}(X^{1/2} (\log X)^{3/2} (\log \log X)^{1/2} \log\log \log{X}).$$ 
We finish the proof of Theorem \ref{pointline} by noting that $B$ must also be of full measure, as it is the intersection of finitely many sets of full measure.

\subsection{The point-cusp problem}
\subsubsection{Proof of Theorem \ref{pointcusp}}
\label{pointcuspsec}
We want to study the asymptotics of
$$
N(X;z,\mathfrak a)
:=
\#\left\{
\gamma\in\Gamma_{\mathfrak a}\backslash\Gamma \mid
\operatorname{Im}(\sigma_{\mathfrak a}^{-1}\gamma z)\ge X^{-1}
\right\}.
$$
We use the convention $N(0;z,\mathfrak{a})=0$.
We consider the spectral expansion of the counting function, $$N(X;z,\mathfrak{a})=M(X;z,\mathfrak{a})+E(X;z,\mathfrak{a}),$$
where $M$ is the contribution of the discrete spectrum, including residual, and $E$ is the contribution of the continuous spectrum.

The contribution of cuspidal Maa\ss{} eigenfunctions in $M$ is zero.  Therefore,
\begin{equation} \label{maincusp}M(X;z,\mathfrak{a})=\frac{X}{\operatorname{vol}(\Gamma\backslash\mathbb H^2)}
+
\sum_{s_j\in\mathcal P_{\mathfrak{a}}}
\frac{\operatorname*{Res}_{s=s_j} \mathcal{E}_{\mathfrak{a}}(z,s)}
{s_j}\,X^{s_j}.\end{equation}
Furthermore, in $L^2(\Gamma \backslash \H)$, we have 

\begin{equation} \label{contprojection}E(X,z):=E(X;z,\mathfrak{a})=\frac{1}{4\pi}\sum_{\mathfrak{b}}\int_{\mathbb{R}}\left\langle N(X;\cdot,\mathfrak{a}),
\mathcal{E}_{\mathfrak b}(\cdot,\tfrac12+it)
\right\rangle \mathcal{E}_{\mathfrak{b}}(z,\tfrac{1}{2}+it) dt. \end{equation}

For $Y \geq X \geq 0$, write

$$
f_{X,Y}(y)=\mathbbm{1}_{[1/Y,\,1/X)}(y),
$$
where $1/0$ is interpreted as $\infty$, and consider the associated incomplete Eisenstein series, defined by
$$F_{\mathfrak{a}}(z):=\sum_{\gamma\in\Gamma_{\mathfrak a}\backslash\Gamma}
f_{X,Y}\bigl(\hbox{Im}(\sigma_{\mathfrak a}^{-1}\gamma\,z)\bigr).$$

It is easy to see that 
$$N(Y;z, \mathfrak{a})-N(X;z, \mathfrak{a})=F_{\mathfrak{a}}(z). $$

Define further
$\Delta_{(X,Y]} E(z)=E(Y,z)-E(X,z).$
By equation (\ref{contprojection}), we have
$$
\Delta_{(X,Y]} E(z)
=\frac{1}{4\pi}\sum_{\mathfrak{b}}\int_{\mathbb{R}}\left\langle F_{\mathfrak{a}},
\mathcal{E}_{\mathfrak b}(\cdot,\tfrac12+it)
\right\rangle \mathcal{E}_{\mathfrak{b}}(z,\tfrac{1}{2}+it) dt. $$
Hence, by Parseval, we have

\begin{equation} \label{parseval}
\|\Delta_{(X,Y]}E\|_2^2
=
\frac1{4\pi}\sum_{\mathfrak b}
\int_{\mathbb R}
\left|
\left\langle F_{\mathfrak{a}},
\mathcal{E}_{\mathfrak b}(\cdot,\tfrac12+it)
\right\rangle
\right|^2\,dt.
\end{equation}

By unfolding, we have 

$$\left\langle F_{\mathfrak{a}},
\mathcal{E}_{\mathfrak b}(\cdot,\tfrac12+it)
\right\rangle = \int_{0}^{+\infty} f_{X,Y}(y)  \int_{0}^{1} \overline{\mathcal{E}_{\mathfrak{b}}\left(\sigma_{\mathfrak{a}}(x+iy),\tfrac{1}{2}+it\right)}dx \frac{dy}{y^2}.$$
Using the scattering matrix convention from \cite[\S 6.3.]{Iwaniec2002}, we can write this as

\begin{align*}\left\langle F_{\mathfrak{a}},
\mathcal{E}_{\mathfrak b}(\cdot,\tfrac12+it)
\right\rangle &= \int_{0}^{+\infty} f_{X,Y}(y)  \left(\delta_{\mathfrak{a}\mathfrak{b}}y^{1/2-it}+\overline{\phi_{\mathfrak{b}\mathfrak{a}}\left(\tfrac{1}{2}+it\right)}y^{1/2+it} \right)\frac{dy}{y^2} \\
&= \delta_{\mathfrak{a}\mathfrak{b}}\hat{f}_{X,Y}\left(\tfrac{1}{2}+it \right)+\overline{\phi_{\mathfrak{b}\mathfrak{a}}\left(\tfrac{1}{2}+it\right)}\hat{f}_{X,Y}\left(\tfrac{1}{2}-it \right),\end{align*}

where

$$\hat{f}_{X,Y}(s):=\int_{0}^{+\infty}f_{X,Y}(y)y^{-s} \frac{dy}{y}$$ is the Mellin transform of $f_{X,Y}$.
Hence, for $s=1/2+it$, we have

$$\sum_{\mathfrak b}
\left|
\left\langle F_{\mathfrak{a}},
\mathcal{E}_{\mathfrak b}(\cdot,s)
\right\rangle
\right|^2 \leq 2\sum_{\mathfrak b}\delta_{\mathfrak{a}\mathfrak{b}}\left|\hat{f}_{X,Y}\left(s \right)\right|^2+2\left(\sum_{\mathfrak b}\left|\phi_{\mathfrak{b}\mathfrak{a}}\left(s\right)\right|^2\right)\left|\hat{f}_{X,Y}\left(1-s \right)\right|^2.$$

Recall now that 
$$\sum_{\mathfrak b}\left|\phi_{\mathfrak{b}\mathfrak{a}}\left(\tfrac{1}{2}+it\right)\right|^2=1,$$
giving

$$\sum_{\mathfrak b}
\left|
\left\langle F_{\mathfrak{a}},
\mathcal{E}_{\mathfrak b}(\cdot,\tfrac12+it)
\right\rangle
\right|^2 \leq 4\left|\hat{f}_{X,Y}\left(\tfrac{1}{2}+it \right)\right|^2.$$
Therefore, by equation (\ref{parseval}), we have

$$
\|\Delta_{(X,Y]}E(\cdot)\|_2^2
 \ll \int_{\mathbb{R}}\left|\hat{f}_{X,Y}\left(\tfrac{1}{2}+it \right) \right|^2 dt.
$$
Using the Mellin--Plancherel theorem \cite[equation~{1.14.38}]{NIST:DLMF}, we conclude that
\begin{equation}\label{pointcuspkey}
\|\Delta_{(X,Y]}E(\cdot)\|_2^2
 \ll 
2\pi\int_0^\infty
|f_{X,Y}(y)|^2\frac{dy}{y^2}
=2\pi(Y-X).
\end{equation}

Now let $n$ be a positive integer. For every integer $m$ with $0 \leq m \leq 2^n$, using the binary expansion of $m$, the interval $(0,m]$ can be partitioned into at most $n+1$ intervals of the form $I(t,j)=(t2^j,t2^j+2^j]$, where $0 \leq j \leq n$ and $0 \leq t < 2^{n-j}$. Write $I_1,\dots, I_{n+1}$ for those intervals, including some empty ones if needed. Then, $\Delta_{(0,m]}E(z)=\sum_{i}\Delta_{I_i} E(z)$. By Cauchy-Schwarz, this gives
$$|\Delta_{(0,m]}E(z)|^2 \leq (n+1) \sum_{i}|\Delta_{I_i} E(z)|^2.$$
Hence,
$$\max \limits_{m \leq 2^n}|\Delta_{(0,m]}E(z)|^2 \leq (n+1) \sum_{j=0}^{n}\sum_{0 \leq t <2^{n-j}}|\Delta_{I(t,j)} E(z)|^2.$$
Using equation (\ref{pointcuspkey}), we get
$$\int_{\Gamma \backslash \H}\max \limits_{m \leq 2^n}|\Delta_{(0,m]}E(z)|^2 d\mu(z)  \ll (n+1) \sum_{j=0}^n2^{n-j} \cdot 2^j=(n+1)^22^n.$$

Using the Chebyshev--Borel--Cantelli argument from the proof of Theorem \ref{mainmaximal}, there exists a set $B$ of full measure, such that, for every $z \in B$ and sufficiently large integer $n$, we have
\begin{equation} \label{cusperrorintegers}\max\limits_{m \leq 2^n}|E(m,z)-E(0,z)| \ll 2^{n/2}n^{3/2}(\log{n})^{1/2}\log\log{n}. \end{equation}

Let $k$ be a nonnegative integer. For $X \in[k,k+1]$,
we have $$N(k;z,\mathfrak{a}) \leq N(X;z,\mathfrak{a}) \leq N(k+1;z,\mathfrak{a}).$$
Furthermore, equation (\ref{maincusp}) gives$$M(X;z,\mathfrak{a})=M(k;z,\mathfrak{a})+O(1)=M(k+1;z,\mathfrak{a})+O(1).$$
Hence,
$$E(k,z)+O(1) \leq E(X,z) \leq E(k+1,z)+O(1).$$
Therefore, using equation (\ref{cusperrorintegers}), we have
$$|E(X,z)| \ll X^{1/2}(\log{X})^{3/2}(\log\log{X})^{1/2}\log\log\log{X},$$
as required.
\subsubsection{Proof of Theorem \ref{omegatheorem}}

Assume that $\Gamma$ is not essentially cuspidal. In other words,
 \begin{equation} \label{essentiallycuspidal} N_{\Gamma}(T):=\#\left\{ j \mid 0\leq \lambda_j \leq T\right\} \not\sim \frac{\operatorname{vol}(\Gamma\backslash\mathbb H^2)}{4 \pi}T. \end{equation}
 
 Recall the scattering-poles form of Weyl law \cite[Theorem 4.13]{muller1992spectral}, stating that 
 \begin{equation} \label{scatteringweyl} N_{\Gamma}(T^2) +N_{p}(T)/2\sim \frac{\operatorname{vol}(\Gamma\backslash\mathbb H^2)}{4 \pi}T^2. \end{equation}
 Here, $N_{p}(T)$ denotes the number of poles $\rho$, with multiplicity, of the determinant of the scattering matrix associated to $\Gamma$, with $|\rho| \leq T$. To avoid ambiguity, we write $(\rho_n)_n$ for the sequence of these poles, with multiplicity.
 Combining equations (\ref{essentiallycuspidal}) and (\ref{scatteringweyl}), we have that
$N_{p}(T)=\Omega(T^2)$. 
On the other hand, Selberg’s weighted pole estimate \cite[Corollary 3.29]{muller1992spectral} (see also \cite[Theorem 2.22]{hejhal1983selberg}) gives us that
$$
\sum_{i: \, \textup{Re}(\rho_i)<1/2}
\frac{\frac12-\hbox{Re}(\rho_i)}{|\rho_i-\frac12|^2}<\infty.
$$
Fix $\epsilon>0$, and let
$$
n_\epsilon(T)
:= \#\bigl\{i \mid |\rho_i|\leq T,\ \hbox{Re}(\rho_i)\leq \tfrac{1}{2}-\epsilon\bigr\}.
$$
For $T\geq R$, each pole counted by $n_\epsilon(T)-n_\epsilon(R)$
contributes at least $\epsilon/(T+\tfrac{1}{2})^2$ to the weighted sum.
Thus,
$$
\frac{n_\epsilon(T)-n_\epsilon(R)}{T^2}
\leq
\frac{(T+\tfrac{1}{2})^2}{\epsilon T^2}
\sum_{\substack{|\rho_i|>R\\ \textup{Re}(\rho_i)<\tfrac{1}{2}}}
\frac{\tfrac{1}{2}-\textup{Re}(\rho_i)}{|\rho_i-\tfrac{1}{2}|^2}.
$$
Taking the limit superior as $T\to\infty$ and then letting $R\to\infty$
gives $n_\epsilon(T)=o(T^2)$. Therefore, $N_p(T)=\Omega(T^2)$ forces a sequence of distinct, nonreal poles $\rho_n'$ with $\hbox{Re}(\rho_n') \uparrow 1/2$ and $|\rho'_n| \rightarrow \infty$.
Recall now that every pole of the scattering determinant must be a pole of an entry, and therefore of an Eisenstein series. As there are finitely many cusps, after passing to a subsequence, we can assume that all the $\rho'_n$ are poles of one fixed Eisenstein series $\mathcal E_{\mathfrak a}(z,s)$. In fact, there is a full measure set $B \subset \Gamma \backslash \H$ such that, for every $z \in B$ and every $n$, the complex number $\rho'_n$ is a pole of $\mathcal E_{\mathfrak a}(z,s)$. This is because the leading Laurent coefficient at $\rho'_n$ is a real analytic function of $z$, and therefore its zero set has measure zero.

We want to show that, for this choice of $\mathfrak{a}$, and for any $z \in B$, we have $E(X;z,\mathfrak a)=\Omega(X^{1/2-\delta})$ for every $\delta>0$. Suppose the contrary. In particular, suppose that there is some $z \in B$ and some $\delta>0$ such that $E(X;z,\mathfrak a)=O(X^{1/2-\delta})$. We note that, by the definition of $(\rho'_n)_n$, there must be some $m$ such that $\hbox{Re}(\rho'_m)>1/2-\delta$. Let $\rho':=\rho'_m$.

Consider now the function
$$E_{z}(s):=\int_1^ \infty E(X;z,\mathfrak a)X^{-s-1}\,dX. $$
On one hand, as $E(X;z,\mathfrak a)=O(X^{1/2-\delta})$, we have that $E_{z}(s)$ defines a holomorphic function throughout $
\hbox{Re}(s)>1/2-\delta.
$
On the other hand, $\rho'$ is a pole of $E_z(s)$. This contradicts the above.
To see why $\rho'$ is a pole of $E_z(s)$, write
\begin{align*}E_{z}(s)&=\int_{0}^{\infty} N(X;z,\mathfrak a)X^{-s-1}\,dX-\int_{0}^{1} N(X;z,\mathfrak a)X^{-s-1}\,dX-\int_{1}^{\infty} M(X;z,\mathfrak a)X^{-s-1}\,dX \\ &= \frac{\mathcal E_{\mathfrak a}(z,s)}s-\int_0^1N(X;z,\mathfrak a)X^{-s-1}\,dX-\frac{(\operatorname{vol}(\Gamma\backslash\mathbb H^2))^{-1}}{s-1}
-\sum_{s_j\in\mathcal P_{\mathfrak{a}}}
\frac{\operatorname*{Res}_{s=s_j} \mathcal{E}_{\mathfrak{a}}(z,s)}
{s_j(s-s_j)}. \end{align*}
For $z \in B$,  we have that $N(X;z,\mathfrak a)=0$ for sufficiently small $X>0$, and hence $$\int_0^1N(X;z,\mathfrak a)X^{-s-1}\,dX$$
is an entire function. In particular, the pole $\rho'$ is preserved. 
Here, we used the fact that, for $\hbox{Re}(s)>1$,
$$
\int_0^\infty N(X;z,\mathfrak a)X^{-s-1}\,dX
=\frac{\mathcal E_{\mathfrak a}(z,s)}s,
$$
and then extended meromorphically.

\section{ Point-hyperplane problem} \label{secth3}
\subsection{Background theory} \label{phpbackground}
Let $\Gamma$ be a discrete cocompact subgroup of $\hbox{PSL}_2(\mathbb{C})$ and $\mathcal{P}$ be a totally geodesic hyperplane in $\Hh $. The totally geodesic hypersurfaces in $\Hh $ are the Euclidean hyperplanes and hemispheres orthogonal to the complex plane.
    By conjugation, we can assume that $\mathcal{P}=\left\{ (x_1,0,y) \mid y>0\right\}$ where, for $p \in \Hh $, we write $p=(z,y)=(x_1,x_2,y),$
where $y>0$. For $\gamma \in \pslc$, we define
$$\gamma p=\left(\frac{(az+b)\overline{(cz+d)}+a\overline{c}y^{2}}{\|cp+d\|^{2}},\frac{y}{\|cp+d\|^{2}}\right).$$ 

 The stabiliser of $\mathcal{P}$ in $\pslc$ is
        $$\hbox{Stab}_{\pslc}(\mathcal{P})=\pslr\bigcup\left(\begin{smallmatrix}i &\\ &-i\end{smallmatrix}\right)\pslr.$$
Equipped with the hyperbolic metric, the hyperbolic distance between
    $p,\,p'$ is given by
    $$\cosh d_{\mathbb{H}^3}(p,p')=\delta(p,p'),$$
    where $\delta$ is the standard point-pair invariant.
    
 As in \cite[\S 2]{niko}, we define the coordinates
    \begin{align*}
        x &= x_{1}, &
        u &= \log \sqrt{x_{2}^{2}+y^{2}}, &
        v &= \arctan\frac{x_{2}}{y}.
    \end{align*}
    We now write $p=(x(p), u(p), v(p))$. With respect to these coordinates, we have

    \begin{equation*}
            ds^{2} = \frac{dx^{2}}{e^{2u}\cos^{2}v}+\frac{du^{2}+dv^{2}}{\cos^{2}v}, \quad
            d\mu(p) = \frac{dx\,du\,dv}{e^{u}\cos^{3}v},
        \end{equation*}
        and
        $$\Delta = e^{2u}\cos^{2}v\pd[2]{x} + \cos^{2}v\left(\pd[2]{u}+\pd[2]{v}\right)-\cos^{2}v\pd{u}+\sin v\cos v\pd{v}.$$

        Furthermore, in the new coordinate system, we have $\mathcal{P}=\left\{ (x,u,0) \right\}.$

        It is easy to see that
        $$\cosh \left(d_{\Hh}(p, \mathcal{P})\right)=(\cos(v(p)))^{-1},$$
and, hence, $$ N(X;p, \mathcal{P})
:=\#\{\gamma\in H\backslash\Gamma \mid \cosh d_{\Hh }(\gamma p,\mathcal P)\leq X\}=\#\{\gamma\in H\backslash\Gamma \mid
(\cos v(\gamma p))^{-1}\leq X\}.
$$

Analogously to equation (\ref{hubers}), for a piecewise continuous function $f$ having compact support in $[1,\infty)$,  we define
    $$A(f)(p)= \sum_{\gamma\in H \backslash \Gamma}f\left(\frac{1}{\cos^{2}v(\gamma p)}\right),$$
    so that, for $f_X(q)$ the indicator of $\left\{q \leq X^2\right\}$, we have $A(f_X)(p)=N(X;p,\mathcal{P})$.
    
    Let $(u_j)_j$ be a complete orthonormal system of real-valued automorphic eigenfunctions of $-\Delta$, with corresponding eigenvalues $(\lambda_j)_j$. We write $\lambda_j=s_j(2-s_j)$ and $s_j=1+it_j$, where $t_j$ is either nonnegative real or purely imaginary. More precisely, either  $\hbox{Re}(s_j)=1$ and $\hbox{Im}(s_j) \geq 0$, or $s_j \in [1,2]$. Define the corresponding periods via 
\begin{equation}\label{nikoperiodsdef}\hat{u}_j=\int_{H\backslash\mathcal{P}}u_{j}(x,u,0)\frac{du\,dx}{e^{u}}. \end{equation}
   In \cite[Lemma 3.1]{niko}, Laaksonen derives the following spectral expansion.
    \begin{lemma}[{\cite[Lemma 3.1]{niko}}]\label{lem:specexp}
      Let $f$ be a compactly supported, continuously differentiable, piecewise $C^2$ real function. We have
    $$A(f)(p) = 2\sum_{j} c(f,t_{j}) \hat{u}_{j}u_{j}(p),$$
         where
        $$c(f,t_{j})=\int_{0}^{\frac{\pi}{2}}f\left(\frac{1}{\cos^{2}v}\right)\frac{\xi_{\lambda_{j}}(v)}{\cos^{3}v}dv,$$
        where $\xi_{\lambda}$ is the solution of the ordinary differential equation
        $$\cos^{2}v\,\xi_{\lambda}''(v)+\sin v\cos v\,\xi_{\lambda}'(v)+\lambda\xi_{\lambda}(v)=0,$$
        with the initial conditions
        $$\xi_{\lambda}(0)=1,\quad\xi_{\lambda}'(0)=0.$$
    \end{lemma}

    Furthermore, he proves the following analogues of Lemmas \ref{lemmaweyl} and \ref{lemmahuber}.
\begin{lemma}[{\cite[Lemma 3.5]{niko}}] \label{weylniko} For $T\geq1$, we have
$$\sum_{0\leq t_j\leq T}|u_j(p)|^2
\ll_\Gamma T^3,
$$
uniformly in $p\in\Gamma\backslash \Hh$.
\end{lemma}
    \begin{lemma}[{\cite[Lemma 3.3]{niko}}]\label{lem:period} 
        For $T \geq 1$, we have
        $$ \sum_{|t_{j}|\leq T}|\hat{u}_{j}|^{2}\ll T.$$
    \end{lemma}

\subsection{Proof of Theorem \ref{pointhyperplane}}

Following Laaksonen \cite[eq. (6)]{niko}, we let
\begin{equation} \label{azp}
M(X,p)=
\frac{\operatorname{vol}(H\backslash\mathcal P)}
{\operatorname{vol}(\Gamma\backslash \Hh )}X^2
+\sum_{1<s_j<2}\frac{2^{s_j-1}}{s_j}
\hat{u}_j u_j(p)X^{s_j},
\end{equation}
and
\begin{equation*}
E(X;p, \mathcal{P})=N(X;p, \mathcal{P})-M(X,p).
\end{equation*}

For $X=e^r$, $r\geq1$, set
$$
U=(X^2-1)^{1/2},\qquad
A=\operatorname{arcsinh}U=\operatorname{arcosh}X,
\qquad
\delta=X^{-1}.
$$
Following \cite[\S 3]{niko}, we construct comparison functions by convolution in the normal coordinate
$$
w=\operatorname{arcsinh}(\tan v),
\qquad
\cosh w=(\cos v)^{-1}.
$$
Define
$$
\chi_\delta(w)=\frac{1}{2\delta}\mathbf1_{[-\delta,\delta]}(w),
\qquad
\overline F_X^\pm(w)
=
\mathbf1_{[-(A\pm2\delta),\,A\pm2\delta]}(w),
$$
and put
$$
F_X^\pm=\overline F_X^\pm*\chi_\delta*\chi_\delta,
\qquad
f_X^\pm(q)=F_X^\pm(\operatorname{arcosh}\sqrt q),
\quad q\geq1.
$$
Here $A>2\delta$, since $r\geq1$. These functions satisfy
$$
0\leq f_X^-\leq\mathbf1_{[1,X^2]}\leq f_X^+\leq1,
$$
and hence
$$
A(f_X^-)(p)\leq N(X;p,\mathcal P)\leq A(f_X^+)(p).
$$
Defining
$$
E^\pm(X,p)=A(f_X^\pm)(p)-M(X,p),
$$
we therefore have
$$
E^-(X,p)\leq E(X;p,\mathcal P)\leq E^+(X,p).
$$
It is enough to prove the required bound for $E^+$ and $E^-$. 

By Lemma \ref{lem:specexp}, we have
$$
A(f_X^\pm)(p)
=
2\sum_j c(f_X^\pm,t_j)\widehat u_j u_j(p).
$$
The convergence needed to use this expansion pointwise will be verified below. We first estimate the transforms and their derivatives.

Write
$$
d(f,t)=\int_{\mathbb R}f(\cosh^2w)\cosh(sw)\,dw, \quad s=1+it.
$$
The transform identity and convolution formula in \cite[Lemma 3.4]{niko} give
$$
4c(f,t)=d(f,t)+d(f,-t),
\qquad s=1+it,
$$
and
\begin{equation} \label{exactniko}
d(f_X^\pm,t)
=
\frac{2\sinh(s(A\pm2\delta))}{s}
\left(\frac{\sinh(s\delta)}{s\delta}\right)^2.
\end{equation}
At $s=0$, the right-hand side is interpreted by continuity.

\begin{lemma} \label{cfxtestim} Let $X=e^r$, $r\geq1$, and let $f_X^\pm$ be as above.
\begin{itemize}
\item[\textup{(a)}] For a fixed real number $0\leq t\leq1$, we have
$$
|c(f_X^\pm,t)|\ll_t X.
$$

\item[\textup{(b)}] For fixed $s=1+it\in(1,2]$, we have
$$
2c(f_X^\pm,t)
=
\frac{2^{s-1}}{s}X^s+O_s(X).
$$

\item[\textup{(c)}] For real $t\geq1$, set
$$
Q(t,r)=\min\{1,(X/t)^{2}\}.
$$
Then, uniformly for $r,t\geq1$,
$$
X^{-1}|c(f_X^\pm,t)|\ll t^{-1}Q(t,r)
$$
and
$$
\left|\partial_r\bigl(X^{-1}c(f_X^\pm,t)\bigr)\right|
\ll Q(t,r).
$$
The derivative includes the dependence $\delta=e^{-r}$.
\end{itemize}
\end{lemma}
\begin{proof} Put
$$
\alpha_\pm=A\pm2\delta,
\qquad
W(s,\delta)=\left(\frac{\sinh(s\delta)}{s\delta}\right)^2.
$$
Part (a) follows directly from equation (\ref{exactniko}).

We note that, for a fixed real number $s$, we have
$$
W(s,\delta)=1+O_s(\delta^2).
$$
On the other hand, for every positive real $s$, we have
$$\frac{2\sinh(s\alpha_{\pm})}{s} =\frac{e^{s\alpha_{\pm}}-e^{-s(A\pm 2\delta)}}{s}=\frac{2^s}{s}X^s+O_s(X^{s-1}+X^{-s}).$$
Furthermore, for $s \rightarrow 0$, we have that 
$$\frac{2\sinh(s\alpha_{\pm})}{s} \rightarrow 2\alpha_{\pm} \ll \log{X}.$$
Hence, combining the above estimates with equation (\ref{exactniko}), we have that, for any fixed real number $1<s \leq 2$,

$$
2c(f_X^\pm,t)
=
\frac{2^{s-1}}{s}X^s
+
O_s\bigl(X^{s-1}+X^{2-s}+\log{X}\bigr)
=
\frac{2^{s-1}}{s}X^s+O_s(X).
$$
This completes the proof of part (b).

For (c), take $s=1+it$, $t\geq1$. Since
$$
|\sinh(s\delta)|^2=\sinh^2\delta+\sin^2(t\delta),
$$
we have
$$
|W(s,\delta)|
\ll \min\{1,(t\delta)^{-2}\}
=
Q(t,r).
$$
Also,
$$
|\sinh(s\alpha_\pm)|+|\cosh(s\alpha_\pm)|
\ll e^A\ll X.
$$
The exact transform formula therefore gives
$$
X^{-1}|d(f_X^\pm,t)|
\ll t^{-1}Q(t,r),
$$
which proves the first estimate in (c).

To estimate the derivative, note that
$$
|\partial_rW(s,\delta)|
\ll\min\{1,(t\delta)^{-1}\}.
$$
Indeed, we have
$$
\partial_rW(s,\delta)=2W(s,\delta)-\frac{\sinh(2s\delta)}{s \delta}
\ll\min\{1,(t\delta)^{-1}\}.
$$

Moreover,
$$
\alpha_\pm'
=
\frac{X}{\sqrt{X^2-1}}\mp2\delta
\ll1.
$$
Differentiating the exact transform formula, we obtain
$$
\begin{aligned}
\partial_r\bigl(X^{-1}d(f_X^\pm,t)\bigr)
={}&-X^{-1}d(f_X^\pm,t)\\
&+2X^{-1}\alpha_\pm'\cosh(s\alpha_\pm)W(s,\delta)\\
&+\frac{2X^{-1}\sinh(s\alpha_\pm)}{s}\,
\partial_rW(s,\delta).
\end{aligned}
$$
Therefore,
$$
\left|\partial_r\bigl(X^{-1}d(f_X^\pm,t)\bigr)\right|
\ll
t^{-1}Q(t,r)+Q(t,r)
+t^{-1}\min\{1,(t\delta)^{-1}\} \ll Q(t,r).
$$
Taking real parts proves the derivative estimate. 
\end{proof}
By Lemma \ref{cfxtestim}(a),(b), we have that
\begin{equation} \label{Etocuu}
E^\pm(X,p)
=
2\sum_{t_j\geq1}c(f_X^\pm,t_j)\widehat u_j u_j(p)+O_{\Gamma,\mathcal P}(X),
\end{equation}
uniformly in $p\in\Gamma\backslash \Hh$. Indeed, there are only finitely many small eigenvalues and real parameters $0\leq t_j<1$, their eigenfunctions are bounded on the compact space $\Gamma\backslash \Hh$, and
$$
\widehat u_0u_0
=
\frac{\operatorname{vol}(H\backslash\mathcal P)}
{\operatorname{vol}(\Gamma\backslash \Hh)}.
$$
The pointwise convergence assertion is justified in the next lemma.

Fix a sign and reindex the spectral parameters $t_j\geq1$. We will apply Theorem \ref{mainmaximal} with
$$
\phi_j=u_j,
\qquad
c_j(r)=2e^{-r}c(f_{X}^\pm,t_j)\widehat u_j,
\qquad
D(t,r)=Q(t,r).
$$
The exact transform formula shows that $c_j\in C^1([1,\infty))$.

\begin{lemma} \label{planeadmissble} For these choices of $t_j,c_j,\phi_j$ and $D$, conditions (C0), (C1), (C2) and (C3) hold. Moreover, for $r$ large enough,
$$
s(r):=\sum_{k\geq0}D(2^k,r)\ll r.
$$
\end{lemma}
\begin{proof} 

Condition (C0) follows from combining Lemma \ref{weylniko}, Lemma \ref{lem:period} and Lemma \ref{cfxtestim}(c), in a similar manner as Lemma \ref{conditioncheck}.

For condition (C1), Lemma \ref{cfxtestim}(c) gives
$$
|t_jc_j(r)|+|\partial_rc_j(r)|
\ll Q(t_j,r)|\widehat u_j|.
$$
Since $Q(t,r)$ is non-increasing in $t$, Lemma \ref{lem:period} yields
$$
\begin{aligned}
\sum_{T\leq t_j<2T}
\left(|t_jc_j(r)|^2+|\partial_rc_j(r)|^2\right)
&\ll
Q(T,r)^2
\sum_{T\leq t_j<2T}|\widehat u_j|^2\\
&\ll_{\Gamma,\mathcal P}TQ(T,r)^2,
\end{aligned}
$$
as required.

For condition (C2), if $a\leq r\leq a+1$, then
$$
D(t,r)
=
\min\{1,e^{2r}t^{-2}\}
\leq
e^2\min\{1,e^{2a}t^{-2}\}
=
e^2D(t,a).
$$

Finally,
$$
\sum_{k\geq0}D(2^k,r)
=
\sum_{2^k\leq e^r}1
+
e^{2r}\sum_{2^k>e^r}2^{-2k}\ll1+r \ll r.
$$
This proves (C3) and the asserted bound for $s(r)$.

\end{proof}
We let
$
g(x)=x^{3/2}(\log x)^{1/2}\log\log x
$ for $x \geq 4$,
and $g(x)=g(4)$ for $1\leq x<4$. Since $s(r)\ll r$, the same summability argument as in equation (\ref{gadmissible}) allows us to apply Theorem 8 for this choice of $g$. Combining this with equation (\ref{Etocuu}), we deduce that there is
a full-measure set
$Z_\pm\subset\Gamma\backslash \Hh$ such that, for every $p\in Z_\pm$,
$$
E^\pm(X,p)
=
O_{\Gamma,\mathcal P,p}\bigl(X+Xg(\log X)\bigr)
\qquad(X\longrightarrow\infty).
$$
Let $Z=Z_+\cap Z_-$. This set has full measure. For $p\in Z$, the comparison
$$
E^-(X,p)\leq E(X;p,\mathcal P)\leq E^+(X,p)
$$

gives
$$
E(X;p,\mathcal P)
=
O_{\Gamma,\mathcal P,p}\!\left(
X(\log X)^{3/2}(\log\log X)^{1/2}\log\log\log X
\right).
$$
This proves Theorem \ref{pointhyperplane}.

\subsection{Proof of Theorem \ref{meanvaluethm}}
Define
\begin{equation} \label{bdefin} e(X; p, \mathcal{P}):=E(X;p,\mathcal{P})-\sum_{1 < s_j < 2} \frac{2^{1-s_j}}{2-s_j}\hat{u}_ju_j(p)X^{2-s_j}-\frac{\operatorname{vol}(H\backslash\mathcal P)}
{\operatorname{vol}(\Gamma\backslash \Hh )}\log{X}, \end{equation}
and $$
V(X, \mathcal{P}):=\int_{\Gamma\backslash \Hh }|e(X;p,\mathcal{P})|^2 d\mu(p).$$
By Parseval, we have that
\begin{equation*} \label{parsevalformv}
V(X,\mathcal{P})
=\sum_j|\langle e(X;\cdot, \mathcal{P}),u_j\rangle|^2.
\end{equation*}
Recall that for $f_X$ the indicator of the interval $[1,X^2]$, we have $A(f_X)(p)=N(X;p,\mathcal{P})$. Therefore, from \cite[Lemma 3.1]{niko}, we have that
\begin{equation} \label{coeffsformv} \langle N(X;\cdot, \mathcal{P}),u_j\rangle=2c(f_X,t_j)\hat{u}_j. \end{equation}
For $X=\cosh{R}$, \cite[Lemma 3.4]{niko} gives
\begin{equation} \label{tformofcft}
c(f_X,t)=
\frac{\bigl(\sinh R\cos(tR)+t\cosh R\sin(tR)\bigr)}
     {1+t^2}.
\end{equation}
Hence, for $t$ real, we have
$$
c(f_X,t)=\frac{X}{\sqrt{1+t^2}}
\cos\bigl(t\hbox{arcosh}{X}-\arctan t\bigr)
+O\!\left(\frac{X^{-1}}{1+t^2}\right).
$$
Therefore,
\begin{align*}
|c(f_X,t)|^2&=
\frac{X^2}{1+t^2}
\cos^2\bigl(t\hbox{arcosh}{X}-\arctan t\bigr)
+O\!\left((1+t^2)^{-3/2}\right).
\end{align*}
Noting that for $\lambda_j \geq 1$ we have $\langle e(X;\cdot, \mathcal{P}),u_j\rangle=\langle N(X;\cdot, \mathcal{P}),u_j\rangle$, equation (\ref{coeffsformv}) now gives
$$\left|\sum_{\lambda_j \geq 1}|\langle e(X;\cdot, \mathcal{P}),u_j\rangle|^2-X^2T(X,\mathcal{P})\right| \ll \sum_{t_j \geq 0}\frac{|\hat{u}_j|^2}{(1+t_j^2)^{3/2}} \ll 1,$$
where the last estimate follows from dyadic partition and Lemma \ref{lem:period}. Therefore,
\begin{equation} \label{vminus} |V(X,\mathcal{P})-X^2T(X,\mathcal{P})|\ll 1+\sum_{ 0 \leq \lambda_j < 1}|\langle e(X;\cdot, \mathcal{P}),u_j\rangle|^2. \end{equation}
Note that, for $s=1+it \in (1,2)$, equation (\ref{tformofcft}) becomes
$$
2c(f_X,t)
=\frac{\sinh(sR)}s+
  \frac{\sinh((2-s)R)}{2-s}=\frac{2^{s-1}}{s}X^{s}+\frac{2^{1-s}}{2-s}X^{2-s}
+O_{s}(1).
$$
Furthermore, for $s=2$, i.e. $t=-i$, we have
$$
2c(f_X,-i)
=\frac{\sinh(2R)}{2}+R=X^2+\log{X}+O(1).
$$
Therefore, we have
\begin{align*}&\sum_{ 0 \leq \lambda_j < 1}|\langle e(X;\cdot, \mathcal{P}),u_j\rangle|^2 \\
&=\left|2c(f_X,-i)-X^2-\log{X}\right|^2|\hat{u}_0|^2+\sum_{1 < s_j <2 }\left|2c(f_X,t_j)-\frac{2^{s_j-1}}{s_j}X^{s_j}-\frac{2^{1-s_j}}{2-s_j}X^{2-s_j}\right|^2|\hat{u}_j|^2 \\
&=\left|O_{\Gamma, \mathcal{P}}(1)\right|^2+\sum_{1 < s_j <2 }\left|O_{s_j}(1)\right|^2|\hat{u}_j|^2
\ll_{\Gamma, \mathcal{P}}1.
\end{align*}
We conclude the proof of Theorem \ref{meanvaluethm} by combining this with equation (\ref{vminus}).

Dividing the formula in Theorem \ref{meanvaluethm} by $X^2$ gives
$X^{-2}V(X,\mathcal{P})=\mathcal{T}(X,\mathcal{P})+o(1)$,
so Corollary \ref{thmmvlimit} reduces to showing that $\mathcal{T}(X,\mathcal{P})$
is bounded and has no limit as $X\to\infty$.
\begin{proof}[Proof of Corollary \ref{thmmvlimit}]
Recall that
    $$T(\cosh{r},\mathcal{P})=
4\sum_{\lambda_j\ge1}
\frac{|\widehat u_j|^2}{1+t_j^2}
\cos^2\!\bigl(t_jr-\arctan t_j\bigr).
$$
Suppose, for contradiction, that $T(X, \mathcal{P}) \rightarrow L$ as $X \rightarrow \infty$.  We note that $$\|T(\cosh{r},\mathcal{P})\|_{\infty} \ll
\sum_{\lambda_j\ge1}
\frac{|\widehat u_j|^2}{1+t_j^2} \ll 1,
$$
and therefore $L$ must be a real number.

There is at least one $j$ with $\lambda_j > 1$ and $\hat{u}_j \neq 0$. Indeed, otherwise the hypersurface measure would have a finite spectral expansion and hence a smooth density, which is impossible. Fix one such $j$ and let $\tau=t_j$. 
We have that
\begin{equation} \label{limit1}
\lim_{R\to\infty}\frac1R
\int_0^R T(\cosh r,P)e^{-2i\tau r}\,dr
=
\frac{e^{-2i\arctan\tau}}{1+\tau^2}
\sum_{j:t_j=\tau}|\widehat u_j|^2\ne0.
\end{equation}
Here, as $R \rightarrow \infty$, the normalized integral extracts the frequency $2\tau$, and all terms at this frequency have the same phase. On the other hand, as $T(X, \mathcal{P}) \rightarrow L$, the standard Ces\`aro averaging argument gives that 
\begin{equation} 
\lim_{R\to\infty}\frac1R
\int_0^R T(\cosh r,P)e^{-2i\tau r}\,dr
=0,
\end{equation}
which contradicts equation (\ref{limit1}).
\end{proof}
\section*{Acknowledgements}
The author was supported by the MTA--HUN-REN RI Lend\"{u}let Automorphic Research Group. The author would like to thank Dimitris Chatzakos for the helpful discussions and his guidance. 

\section*{Statement on the use of Artificial Intelligence}
The author used OpenAI’s ChatGPT (GPT-6 Astra) to discuss and critically review portions of the mathematical
arguments, check calculations, assist with locating references, and improve the exposition. The author independently verified all the calculations, proof steps and references included in the final manuscript and takes full responsibility for its contents.

\renewcommand*{\bibfont}{\small}
\printbibliography
\end{document}